\documentclass[10pt, openright]{article}
\usepackage{amsfonts, amsmath, amssymb, graphics, color}
\usepackage{arev}
\usepackage[a4paper]{geometry}
\usepackage{exact-computation-over-topological-spaces-arxiv}
\usepackage{latexsym}
\usepackage{graphicx}
\usepackage{url}
\usepackage{adjustbox}
\usepackage{authblk}
\usepackage{fancyhdr}
\usepackage{enumitem}
\usepackage{footmisc}
\makeatletter
\renewcommand{\@pnumwidth}{2.25em}
\renewcommand{\@tocrmarg}{3.25em}
\makeatother

\begin{document}
\title{Exact computation over topological spaces: constructive theory and practice}
\date{\normalsize{December 7, 2013}}
\author{Frank Waaldijk\thanks{www.fwaaldijk.nl/mathematics.html\\
\emph{2012 ACM CCS:} [\bol{Mathematics of computing}]: Continuous mathematics--Topology--Point-set topology}}
\maketitle

\small\begin{center}\bol{Abstract}\end{center}
\parr
We give an exposition of Natural Topology (\NTop), which highlights its advantages for exact computation. The \NTop-definition of the real numbers (and continuous real functions) matches the recommendations for exact real computation in \BauKavto\ and \BauKav. We derive similar and new results on the efficient representation of continuous real-valued functions (defined on a suitable topological space). This can be generalized to continuous functions between suitable topological spaces. Other than in \BauKav, we do not need Markov's Principle (but for practice this is a cosmetic difference). 

\parr

\NTop\ is a conceptually simple theory of topological spaces in \bish. It combines a pointfree with a pointwise approach, by integrating the partial-order on the topological basis with a pre-apartness relation. Simpler than formal topology and (constructive) domain theory, \NTop\ enables a smooth transition from theory to practice. We define `natural reals' as sequences of `sufficiently shrinking' rational intervals. The construction directly yields an apartness topology which is equivalent to the metric topology. 

\parr

A similar construction works for all `effective' quotient spaces of Baire space (obtained through a $\Sigma_0^1$-apartness). We work with a countable set of `basic dots' which are usually basic neighborhoods. This allows for an efficient representation of compact spaces by finitely branching trees (contrasting to the framework of formal topology). For the natural reals, we can suffice with lean dyadic intervals (in $\{[\frac{n}{2^{m}}, \frac{n+2}{2^{m}}]\midd n\inn\Z, \minn\}$). All our spaces are (also) pointwise topological spaces, enabling the familiar pointwise style of \bish\ and \class, as well as an incorporation of earlier constructive work in analysis.

\parr
The concept of `refinement morphism' is seen to adequately capture the notion of `continuous function'. A refinement morphism simply sends basic dots to basic dots, in such a way that `points go to points'. In the case of the reals, every \bish-continuous real function can thus be represented by a morphism sending lean dyadic intervals to lean dyadic intervals.

\parr
For a large class of spaces we prove that continuous functions can be represented by morphisms, that is if we work in \class, \intu\ or \russ. This then should be enough validation also for \bish. We conclude that \NTop\ addresses the need expressed in \BauKavto, to have a framework for constructive topology which is both theoretically and computationally adequate. 

\normalsize
\sectionnb{Theory and practice in constructive topology}\label{introcomp}

\sbsc{A suitable framework for constructive topology}\label{compintro1}
In recent decades, the problem of finding a suitable framework for constructive topology has received increasing attention. Several approaches have been proposed and partly developed, but no definite framework has so far emerged. As an incomplete list, let us mention Domain Theory, Formal Topology, Intuitionistic Topology, Abstract Stone Duality, Apartness Spaces (in the style of Bridges), Type-Two Effectivity and Natural Topology (\WaaNat).

\parr
In \BauKavto, the need is expressed for a unified framework for constructive topology which is both theoretically precise and at the same time suited for exact computation. We think this need for a unified framework is accommodated by Natural Topology (\NTop).

\parr
\NTop\ is developed in \WaaNat, providing a simple constructive framework for topology. \NTop\ is shown to suffice for a \bish\ theory of (separable $T_1$) topological spaces, while accommodating computational practice at the same time. Some simplified computation-related examples are given, but a large part of \WaaNat\ is devoted to proving that the theoretical framework is valid. Also the connections with \class, \intu, \russ\ and formal topology are studied. Especially relevant for the representation of continuous functions is the Lindel\"of-type axiom \BDD\ (defined in \WaaArt) which holds in \class, \intu\ and \russ. Another important theoretical element concerns the development of a simple transfinite inductive machinery, which enables one to work with Heine-Borel properties of compact spaces also in \bish.

\parr
These theoretical issues do not make for easy reading, since the foundations of constructive mathematics are involved in an essential way. It seems worthwhile to accentuate the practical advantages of \NTop\ in a separate paper. We also derive some new results.

\sbsc{Aim and scope of this article}\label{compintro2}
We aim to give a short exposition of Natural Topology (\NTop), which highlights its advantages for exact computation. The proofs of the basic theorems can be found in \WaaNat. In particular we focus on the real numbers. The \NTop-definition of the real numbers (and continuous real functions) matches the recommendations for exact real computation in \BauKavto\ and \BauKav. We derive similar and new results on the efficient representation of continuous real-valued functions (defined on a suitable topological space). This can be generalized to continuous functions between suitable topological spaces. Other than in \BauKav, we do not need Markov's Principle. But for practice this is a cosmetic difference, and even for theory the difference seems inessential. The similarities are more important. We think \NTop\ provides an elegant theoretical frame for the implementations discussed in \BauKav. 

\parr

\NTop\ is a conceptually simple theory of topological spaces in \bish. It combines a pointfree with a pointwise approach, by integrating the partial-order on the topological basis with a pre-apartness relation. Simpler than formal topology and domain theory, \NTop\ enables a smooth transition from theory to practice. We define `natural reals' as sequences of `sufficiently shrinking' rational intervals. The construction directly yields an apartness topology which is equivalent to the metric topology. 

\parr

A similar construction works for all `effective' quotient spaces of Baire space (obtained through a $\Sigma_0^1$-apartness). We work with a countable set of `basic dots' which are usually basic neighborhoods. This allows for an efficient representation of compact spaces by finitely branching trees (contrasting to the framework of formal topology). For the natural reals, we can suffice with lean dyadic intervals (in $\{[\frac{n}{2^{m}}, \frac{n+2}{2^{m}}]\midd n\inn\Z, \minn\}$). All our spaces are (also) pointwise topological spaces, enabling the familiar pointwise style of \bish\ and \class, as well as an incorporation of earlier constructive work in analysis.

\parr
The concept of `refinement morphism' is seen to adequately capture the notion of `continuous function'. A refinement morphism simply sends basic dots to basic dots, in such a way that `points go to points'. In the case of the reals, every \bish-continuous real function can thus be represented by a morphism sending lean dyadic intervals to lean dyadic intervals (this can be generalized to $\R^n$). Using \BDD, we can prove the same for all continuous real functions (which is relevant for \russ). Sometimes, mostly for theoretical purposes, so-called `trail morphisms' are also necessary. Trail morphisms correspond to intuitionistic spread-functions, and in general are less efficient computationally than refinement morphisms.

\parr
For a large class of spaces we can prove that continuous functions can be represented by morphisms, using the axiom \BDD\ (valid for \class, \intu\ and \russ). This then should be enough validation also for \bish. We conclude that \NTop\ addresses the need expressed in \BauKavto, to have a framework for constructive topology which is both theoretically and computationally adequate.

\sbsc{Efficient computation equals efficient representation}\label{compintro3}
We believe that efficient exact computation for topological spaces depends on two key issues, which are closely related. The most important issue seems how to efficiently represent a given topological space. The second issue then becomes how to efficiently represent a continuous function between two efficiently represented spaces. These issues are studied in \WaaNat, also for theoretical purposes. We present the computationally relevant results in this paper.

\sbsc{Structure of this article: theory and practice}\label{compintro4}
The article has two equally important tiers: theory and computational practice. We describe a simple constructive framework for general topology (theory). It takes up some time to illustrate that this framework is comprehensive, notwithstanding its simplicity. We then exploit the simplicity for computational purposes (practice). We repeat some examples given in \WaaNat, and we also add some new theorems.

\parr

By practice we mean: implementations of exact computation. In \BauKavto, \BauKav\ and \KreSpi\ actual implementations of exact real computation are carried out in \pname{RZ}, \pname{HASKELL} and \pname{Coq}. We did not carry out similar implementations for \NTop, but for our real-function morphisms there is a precise correspondence to \BauKavto\ and \BauKav. We agree with the recommendations in \BauKav, one can see \NTop\ as a framework in \bish\ to accomodate these recommendations (conforming to the need expressed in \BauKavto).

\parr
Other than implementations of exact (real) computation, there is applied mathematics to consider. Here the author's knowledge is even less than in regard to exact real computation, so please read the next disclaimer.

\sbsc{Disclaimer}\label{disclaimer}
Applied programming practice is far beyond the author's knowledge. This means that we partly rely on \BauKav\ for our claims of `efficiency' regarding implementations of exact computation. Where applied mathematics in general is concerned, we do not claim any wisdom at all. We think that \NTop\ should be of interest for applied mathematics and applied computation as well, at least as theoretical benchmark, for conceptual reasons, and perhaps also for validation purposes. We have included the (simplified) example of the decision-support system \hwki\ as a conceptual illustration.

\parr
When getting down to the real nitty-gritty, programmers' inventivity and expertise probably lead to other, more efficient solutions. Another reason why our framework for exact computation may not be the most efficient for applied math is that our world in practice is finite. `Infinite-precision arithmetic' actually means `potentially-infinite-precision arithmetic' (the first name is an understandable choice however...). So in practice we always work with finite precision and `rounding' errors. The handling of finite precision can perhaps be done more efficiently when disregarding partial-order properties which we use for potentially-infinite precision, for instance.

\parr
Finally, the volume of research in constructive mathematics is such that we are likely to be unaware of many results and ideas which pertain to our narrative. Ideas in topology like planets often revolve around the same star. Already in \Freu\ one finds an intuitionistic pointfree development, and there is a large body of literature on the subject. Therefore we do not claim wisdom, and certainly not originality, and errors will be gladly corrected when pointed out.

\sectionn{Basic definitions and the natural reals}\label{natbasdef}

\sbsc{Topology first, points later}\label{mathintro1}
The key point of Natural Topology is to define the real numbers -- more generally a (separable $T_1$) topological space -- by \emph{starting}\ with the topology, and obtaining the points of the topological space in the process. Specifically, we start with a countable set of basic dots, which often represent basic neighborhoods in the topology of the space to be constructed.

\sbsc{Pre-natural spaces}\label{nattopdef1}
Our basic mathematical setting involves a countable set \Vt\ of \deff{basic dots}\ of a \deff{natural topological space}\ \Vnatt\ which we build with a number of definitions in this section. Along with the definitions we give some explanations and examples.

\defi
A pre-natural space is a triple \Vprenatt\  where \Vt\ is a countable\footnote{A set $S$ is \deff{countable} iff there is a bijection from \Nt\ to $S$, and \deff{enumerable} iff there is a surjection from \Nt\ to $S$.}\  set of \deff{basic dots}\ and \aprt\ and \leqct\ are binary relations on \Vt, satisfying the properties following below. Here \aprt\ is a \deff{pre-apartness}\ relation (expressing that two dots lie apart) and \leqct\ is a \deff{refinement}\ relation (expressing that one dot is a refinement of the other, and therefore contained in the other). 

\be

\item[(i)] The relations \aprt\ and \leqct\ are \deff{decidable}\ on the basic dots.
\item[(ii)] For all $a, b \inn V$: $a \aprt b$ (`\aat\ is apart from \bbt') if and only if $b \aprt a$. Pre-apartness is symmetric.
\item[(iii)] For all $a \inn V$: $\neg (a \aprt a)$.  Pre-apartness is antireflexive.
\item[(iv)] For all $a, b, c \inn V$: if $a \leqc b$ (`\aat\ refines \bbt') then $c \aprt b$ implies $c \aprt a$. Pre-apartness is \geqc-monotone.
\item[(v)] The relation \leqct\ is a partial order, so for all $ a, b, c\inn V$: $a\leqc a$ and if $a\leqc b\leqc c$ then $a\leqc c$, and if $a\leqc b \leqc a$ then $a=b$. Refinement is reflexive, transitive and antisymmetric. 

\ee
For basic dots we write $a\touch b$ (`$a$ touches $b$') iff $\neg(a\aprt b)$. Then $\touch$ is the decidable complement of $\aprt$.
\edefi

\rem For the motivating example of the real numbers, the basic dots can be thought of as the (closed) rational intervals.\footnote{Using open intervals is also possible, but in general less efficient.}\ Two rational intervals $[a,b]$ and $[c,d]$ are said to be \deff{apart}, notation $[a,b]\aprt[c,d]$, iff either $d\smlr a$ or $b\smlr c$. $[a,b]$ \deff{refines} $[c,d]$, notation $[a,b]\leqc[c,d]$, iff $c\leqq a$ and $b\leqq d$.
\erem

\sbsc{Points arise from shrinking sequences}\label{defpoints}
We turn to infinite sequences (of dots), in order to arrive at points. Looking at our example of rational intervals we see that we need to impose a `sufficient shrinking' condition, otherwise the infinite intersection may contain a whole interval rather than just a point. For an infinite shrinking sequence $\alpha=r_0, r_1, \ldots $ of closed rational intervals ($r_{m+1} \leqc r_{m}$ for all indices $m$) to represent a real number, \alfat\ must `choose' between each pair of apart rational intervals $[a,b]\aprt[c,d]$. By which we mean: for each such pair $[a,b]\aprt[c,d]$, there is an index $m$ such that $r_m \aprt [a,b]$ or $r_m \aprt [c,d]$. (We leave it to the reader to verify that this is indeed equivalent to saying that the infinite intersection of $(r_m)_{\minn}$ contains just one real number.). 

\parr 

The elegance of this approach is that for an infinite shrinking sequence of dots, the property of `being a point' can be expressed by an enumerable condition of pre-apartness. There is no need to talk of `convergence rate' or `Cauchy-sequence', which both presuppose some metric concept. 
To define points, we simply study the real numbers and transfer certain of their nice properties to our general setting.


\defi
A \deff{point}$\!$\ on the pre-natural space \Vprenatt\ is an infinite sequence $p\iz p_0,$ $p_1, p_2 \ldots $ of elements of $V$ that satisfies:

\be
\item[(i)] for all indices $n$ we have: $p_{n+1} \leqc p_{n}$ and there is an index $m$ with $p_m\precc p_n$.
\item[(ii)] If $a, b \inn V$ and $a \aprt b$ then there is an index $m$ such
that $p_m \aprt a$ or $p_m \aprt b$.
\ee
 
Note that any infinite subsequence of $p$ is itself a point (equivalent to $p$ in the natural sense to be defined).  The set of all points on \Vprenatt\  is denoted by \Vcalt.  
\edefi

Since points are infinite sequences, the set \Vcalt\ is generally not enumerable (but all points in \Vcalt\ could be equivalent).

\sbsc{Apartness on points}\label{defaprtpoints}
The points of our pre-natural space \Vprenatt\ are defined, but clearly we obtain many points which are in some sense equivalent (see our example of rational intervals). The constructive approach to an equivalence relation is to look at its strong opposite, namely an \deff{apartness}. 

Therefore it is convenient to extend \aprt\ to points in \Vcalt, and also define when points `belong' to dots, in the obvious way:

\defi
For $p=p_0, p_1,\ldots,  q=q_0, q_1,\ldots \inn \Vcal $ and $a \inn V$:

\be
\item[(i)] $a \aprt p$ and $p \aprt a$ iff $a \aprt p_m$ for some index $m$.
\item[(ii)] $p \aprt q$ iff $p_n \aprt q_n$ for some index $n$.
\item[(iii)] $p\equivv q$ iff $\neg(p\aprt q)$.
\item[(iv)] $p \precc a$ iff $p_m\precc a$ for some index $m$.  This relation is
also referred to as `$a$ is a beginning of $p$' or `$p$ begins with $a$' or `$p$ belongs to $a$'. 
\item[(v)] We write \ahatt\ for the set of all points $p$ such that $a$ is a beginning of $p$. Notice that \ahatt\ is not necessarily closed under \equivv. We write \aclost\ for the \equivv-closure of \ahat. 
\ee
\edefi

It is easy to see that \aprt\ is indeed an apartness (\WaaNat, 1.0.4). In terms of complexity, \aprt\ is a $\Sigma_0^1$-property, whereas \equivvt\ is a $\Pi_0^1$-property. This reflects that apartness of two sequences of dots can be seen at some finite stage, but equivalence of two such sequences is an infinite property. Therefore apartness is better suited for constructive and computational purposes.

\sbsc{Apartness topology is the natural topology}\label{defnattop}
There is a natural topology on the set of points \Vcalt\ of a pre-natural
space \Vprenat. This topology is expressed in terms of apartness and refinement, we call it the \deff{natural topology}\ and also the \deff{apartness topology}, denoted as \Topaprt. \Topaprt\ is the collection of \aprt-open subsets of \Vcalt\ where \aprt-open is defined thus: 

\defi
A set $U\subseteqq \Vcal $ is \deff{\aprt-open} iff for each $x \inn U$ and each $y \inn\Vcal$ we can determine at least one of the following two conditions (they need not be mutually exclusive):

\be
\item[(1)] $y \aprt x$
\item[(2)] there is an index $m$ such that $\hattr{y_m}=\{z\inn\Vcal\midd z \precc y_m\}$ is contained in $U$.
\ee
When the context is clear we simply say `open' instead of `\aprt-open'. 
\edefi

It follows from this definition that an open set is saturated for the equivalence on points (meaning if $U$ is open, $x\inn U$ and $x\equivv y$ then $y\inn U$). We leave this to the reader for easy verification. (This also means that we could replace $\hattr{y_m}$ with $\closr{y_m}$ in (2) above, but in practice this leads to slightly more elaborate proofs). One easily shows that the above indeed defines a topology on \Vcalt\ (else consult \WaaNat, 1.0.5).

\sbsc{Natural spaces}\label{defnatspace}
All the ingredients for our main definition have been prepared. Notice that we did not yet stipulate that each dot should at least contain a point. Also it turns out to be necessary to have a maximal dot, which contains the entire space. These then become the final requirements:

\defi
Let \Vprenat\ be a pre-natural space, with corresponding set of points \Vcalt\ and apartness topology \Topaprt. An element $d$ of $V$ is called a \deff{maximal dot} iff $a\leqc d$ for all $a\in V$. Notice that $V$ has at most one maximal dot\footnote{Actually, it also makes sense to reverse the \leqct-notation, and to consider our maximal dot as being the minimal element, which carries the least information. Then each refinement is `larger' because it carries more information than its predecessor.}, which if existent is denoted \maxdotVt\ or simply \maxdott.  \Vnatt\ is a \deff{natural space}\ iff $V$ has a maximal dot and every $a \inn V$ contains a point.
\edefi

\lem
Let \Vnatt\ be a natural space, derived from the pre-natural space \Vprenat. Let $a\inn V$ be a basic dot. Then the set $\aprt(\hattr{a})=\{z\inn\Vcal\midd z\aprt a\}$ is open in the natural topology. 
\elem

\crl 
For  $x$ in \Vcalt, the set $\{w\inn\Vcal\midd w\aprt x\}$ is open in the natural topology. So a set containing one point (up to equivalence) is closed, showing that every natural space is T$_1$.
\ecrl

\prf see \WaaNat, 1.0.6
\eprf

\rem
The lemma does not imply that basic dots correspond to `closed' subsets in the topology. We can also construct the real numbers as a natural space where the basic dots correspond to open intervals, see remark \ref{defrealnat}. Different representations of the real numbers are relevant for studying computational practice. For example, both decimal reals and `nested-intervals' reals can be represented in \NTop, yielding interesting comparison (see section \ref{appmath}).
\erem

\sbsc{The natural real numbers}\label{defrealnat}
After using the rational intervals as a running example for $V$, we can now formally define the \deff{natural real numbers}\ \Rnatt\ as follows:

\defi
Let $\Rrat\isdef\{[p,q]| p,q \inn \Q | p<q\}\cupp\{(-\infty,\infty)\}$. For two rational intervals $[a,b]$ and $[c,d]$ put $[a,b]\aprtR[c,d]$ iff ($d\smlr a$ or $b\smlr c$) and put $[a,b]\leqcR[c,d]$ iff ($c\leq a$ and $b\leq d$). The maximal dot \maxdotRt\ is obviously $(-\infty,\infty)$. The points on the pre-natural space \Rprenatt\ are called the \deff{natural real numbers}\ (also `natural reals'), the set of natural reals is denoted by \Rnatt. The corresponding natural topology is denoted by \TopRaprt. (Also see the remark later in this paragraph).
\parr
Next, let $\zorrat\isdef \{[p,q]| p,q \inn \Q | 0\leqq p\smlr q\leqq 1\}$, then \zorprenatt\ is a pre-natural space with corresponding natural space $(\zornat, \TopRaprt)$ and maximal dot $\maxdot_{\zor}\iz\zor$. 
\edefi 

\thm
\Rnattopt\ is a natural space which is homeomorphic to the topological space of the real numbers \Rt\ equipped with the usual metric topology.
\ethm  

\prf
Not difficult, see \WaaNat, A.3.0. Notice that by `homeomorphism' we mean the usual definition (a continuous function from one space to the other which has a continuous inverse; `continuous' meaning that the inverse image of an open set is itself open). Also, we are a bit free here, since for a classical theorist we should first move to the quotient space of equivalence classes.
\eprf

\rem Another interesting representation of \Rt\ as a natural space is obtained by changing just very little in the definition. For two intervals $[a,b]$ and $[c,d]$ in $\Rrat$ put $[a,b]\aprtRo[c,d]$ iff ($d\leqq a$ or $b\leqq c$) and put $[a,b]\leqcRo[c,d]$ iff ($c\smlr a$ and $b\smlr d$). Then\Rprenatot\ is a pre-natural space, and the corresponding natural space is again homeomorphic to \Rt. But one sees that the basic dots $[a, b]$ now correspond to the \emph{open} real intervals $(a, b)$.\,\footnote{It would therefore be better to denote the basic dots as open rational intervals $(a, b)$ under this definition.}\ \Rprenatot\ resembles the definition of the formal reals in formal topology (we believe). However, we think compactness is less wieldy in \Rprenatot, which is one reason to stick with \Rprenatt.\footnote{The two spaces are isomorphic in a sense yet to be defined.}
\erem

\sectionn{Natural morphisms}\label{natmorf}

\sbsc{Structure-preserving mappings and real-number representations}\label{whymorphisms}
Insight into natural spaces often comes from mappings which are structure-preserving to some extent. We define different types of such mappings, calling all of them \emph{natural morphisms}. Each natural morphism defines a continuous function with respect to the natural topology. In \class, the structure of natural morphisms sometimes gives a finer distinction between natural spaces, than the structure of continuous functions between their corresponding topological spaces. Then for a classical mathematician, the natural morphisms form an interesting subclass of the class of continuous functions. Still, in \ref{basicneighbor}\ we show that for many natural spaces continuous functions can be represented by a natural morphism (in \class, \intu\ and \russ).

\parr

It turns out there is no \emph{isomorphism} between the natural real numbers and the `natural decimal real numbers', whereas classically these spaces are topologically identical. The natural decimal reals are not closed under basic arithmetic operations, and turn out to be `\pathnat\ connected' but not `\arcnat\ connected'. Similarly, many intuitionistic results can be translated to \NTop, providing an alternative classical way to view important parts of intuitionism.

\parr
Morphisms are suited for efficient computation of continuous functions between topological spaces. The efficiency depends on two factors: the efficient representation (as a natural space) of the spaces involved, and the efficient representation of the continuous function by a refinement morphism (to be defined). 

\sbsc{Different representations of the `same' space}\label{diffrep}
In topology, a central role is played by homeomorphisms. When two spaces are homeomorphic, one can see them as two different representations of the `same' topological space. Yet there is often an intrinsic interest in these different representations. Consider for example \Rt\ and \Rplust. These are two homeomorphic spaces ($(\R, +)$ and $(\Rplus, \cdot)$ are even isomorphic topological groups), but we often have use for one or the other representation, depending on context.

\parr
  
To build an elegant theory and prove its correctness, we look at many different representations of `same' natural spaces. In \NTop\ `sameness' is induced by a special class of natural morphisms called `isomorphisms'. Every isomorphism induces a homeomorphism, but the converse is not true in \class\ (see the above example of the natural decimal real numbers).

\parr

We will give a computationally efficient (isomorphic) representation of the real numbers, which enables continuous functions to be represented by computationally efficient morphisms. We study similar representations for other topological spaces. It seems advantageous to have various representations of a space, which can be used depending on the context. The important thing to note is that \NTop\ gives a unified framework to move from one such representation to another, even efficiently, using isomorphisms.

\sbsc{Natural morphisms 1: refinement morphisms}\label{naturalmorphisms}
We distinguish two types of natural morphisms: \emph{refinement morphisms} (denoted \leqc-morphisms) and \emph{trail morphisms} (denoted \pthh-morphisms). The definition of refinement morphisms is slightly modified from \WaaNat, to fix a minor oversight.

\parr

When going from one natural space to another, a refinement morphism sends basic dots to basic dots, respecting the apartness and refinement relations, in such a way that `points go to points'. This means that any \leqc-morphism\ is an order morphism with respect to the partial order \leqct.\footnote{Not all order morphisms are refinement morphisms though. Our notation ` \leqc-morphism'  can be slightly misleading in this respect.}\ 

\defi 
Let \Vnatot\ and \Wnattot\ be two natural spaces, with corresponding pre-natural spaces \Vprenatot\ and \Wprenattot. Let $f$ be a function from $V$ to $W$. Then $f$ is called a \deff{refinement morphism} (notation: \leqc-morphism) from \Vnatot\ to \Wnattot\ iff for all $a,b\in V$ and all $p=p_0,p_1,\ldots,\ q=q_0,q_1,\ldots\in\Vcal$:

\be
\item[(i)]$f(p)\isdef\ f(p_0), f(p_1), \ldots $ is in \Wcalt\ (`points go to points').
\item[(ii)]$f(p)\aprtto f(q)$ implies $p\aprto q$.
\item[(iii)]$a\leqco b$ implies $f(a)\leqcto f(b)$ (this is an immediate consequence of (i)).
\ee

As indicated in (i) above we write $f$ also for the induced function from \Vcalt\ to \Wcalt. The reader may check that (iii) follows from (i). By (i), a \leqc-morphism $f$ from \Vnatot\ to \Wnattot\ respects the apartness/equivalence relations on points, but not necessarily on dots since $f(a)\aprtto f(b)$ does not necessarily imply $a\aprto b$ for $a,b\in V$. This stronger condition however in practice obtains very frequently.
\edefi

\thm
Let $f$ be a \leqc-morphism from \Vnatot\ to \Wnattot. Then $f$ is continuous.
\ethm

\prf
Easy, see \WaaNat, 1.1.2
\eprf

\sbsc{Refinement morphisms are computationally efficient}
Refinement morphisms are simple in concept. They have the added advantage of computational efficiency. With a suited `lean' representation \sigRt\ of the natural real numbers, \leqc-morphisms from \sigRt\ to \sigRt\ resemble interval arithmetic, and match the recommendations in \BauKav\ for efficient exact real arithmetic.

\parr

More generally, a continuous function between two `lean' natural spaces can usually be represented by a \leqc-morphism (see thm.\,\ref{basicneighbor}, prp.\,\ref{refvstrail2}) and thm.\,\ref{repcontreal}). Therefore we construct `lean' representations of natural spaces (called `spraids'). We believe in the general efficiency of combining \leqc-morphisms with spraids. In this combination, points and continuous functions are of similar type: a sequence of basic dots.

\parr
Spraids turn out to be fundamental for the theroretical perspective as well. To understand the complexities and to prove our framework correct, one needs to define trail morphisms and trail spaces. This can all be found in \WaaNat, we will skip most of the details here.

\sbsc{Natural morphisms 2: trail morphisms}\label{pathmorphisms}
For the most general theoretical setting of natural spaces and pointwise topology, \leqc-morphisms turn out to be too restrictive. This explains our use for the more involved concept of `trail morphism' (denoted \pthh-morphism), defined in the appendix \ref{deftrailmorph}. Trail morphisms play a necessary role in establishing nice properties of natural spaces. Once these properties have been established, we can primarily use \leqc-morphisms\ (see the previous paragraph). Where \leqc-morphisms are defined naturally on basic dots, one can see \pthh-morphisms as mappings which are naturally defined on points.

\parr
The most important property of \pthh-morphisms is that they are technically also \leqc-morphisms, defined on the so-called `trail space' of a natural space (which in turn is a natural space). They are therefore continuous. In this paper we concentrate on refinement morphisms, the interested reader may consult \WaaNat\ for more details on trail morphisms.

\sbsc{Natural morphisms' convention}
The difference between \leqc-morphisms and \pthh-morphisms is often not relevant, which justifies the following:

\conv
If \Vnatot\ and \Wnattot\ are two natural spaces, and $f$ is a \leqc-morphism or a \pthh-morphism from \Vnatot\ to \Wnattot, where the difference is irrelevant, then we simply say: $f$ is a \deff{natural morphism}\ from \Vnatot\ to \Wnattot, or even more simply: a morphism from \Vnatot\ to \Wnattot. 
Only when the difference is relevant will we specify `refinement morphism' and/ or `trail morphism'. This happens mostly in technical proofs or in the context of computation, since refinement morphisms are generally more efficient.
\econv

\sbsc{Composition of natural morphisms}\label{morfcomp}
Given two \leqc-morphisms $f,g$ from natural spaces \Vcalt\ to \Wcalt\ and \Wcalt\ to \Zcalt\ respectively, to form their composition is unproblematic. We leave it to the reader to verify that putting $h(a)\iz g(f(a))$ for all $a\inn V$ defines a \leqc-morphism $h$ from \Vcalt\ to \Zcalt. Composition involving \pthh-morphisms is detailed in \WaaNat, 1.1.6.

\sbsc{Isomorphisms}
We can now define a natural parallel to the topological idea of `homeomorphism'. We will call this parallel `isomorphism'. Isomorphisms between natural spaces will automatically be homeomorphisms, but classically we can find homeomorphic natural spaces which are non-isomorphic. This shows that our theory enriches \class\ as well.

\defi
Let \Vnatot\ and \Wnattot\ be two natural spaces. A natural morphism $f$ from \Vnatot\ to \Wnattot\ is called an \deff{isomorphism}\ iff there is a morphism $g$ from \Wnattot\ to \Vnatot\ such that $g(f(x))\equivvo x$ for all $x$ in \Vcalt\ and $f(g(y))\equivvto y$ for all $y$ in \Wcalt. An isomorphism $f$ from \Vnatot\ to \Vnatot\ is called an \deff{automorphism}\ of \Vnatot, and an \deff{identical automorphism}\ iff $f(x)\equivvo x$ for every $x\inn\Vcal$.
\edefi 

\parr

To see whether certain properties of natural spaces are truly `natural', we check if they are preserved under isomorphisms.

\sectionn{Fundamental natural spaces}

\sbsc{Baire space and Cantor space}
Baire space (\NN) is fundamental because it is a universal natural space (meaning that every natural space can be thought of as a quotient space of Baire space). For Baire space, the relevant partial order \leqct\ is a tree. In \WaaNat, this is exploited to simplify the theory considerably. Cantor space (\zoN) is likewise a universal `fan' by which we mean a space generated by a partial order \leqct\ which is a finitely branching tree. Cantor space can be seen as a universal compact space. 

\sbsc{The class of natural spaces is large}\label{natlarge}
Many spaces can be represented by a natural space. In other words, the class of natural spaces is large. A non-exhaustive and also repetitive list of spaces which can be represented as a natural space:
\be
\item[$\bullet$] every complete separable metric space
\item[$\bullet$] the (in)finite product of natural spaces
\item[$\bullet$] \Nt, \Rt, \Ct, the complex p-adic numbers \Cpt, \RNt, Baire space, Cantor space, Hilbert space \Hit, every Banach space, the space of locally uniformly continuous functions from \Rt\ to \Rt, many other continuous-function spaces, and Silva spaces (see \WaaNat, chapter four).
\ee

Sometimes, classically defined non-separable spaces (for instance function spaces equipped with the sup-norm, see \WaaThe) can be constructed under a different metric to become separable. Although the topology is then not equivalent, one can still work with the space constructively as well. For this, one sometimes needs to construct a completion first, to refind the original space as a subset of the completion. Thinking things through, we do not really see a constructive way to define `workable' spaces other than by going through some enumerably converging process. In this sense we concur with Brouwer. Brouwer's definition of spreads in essence parallels  the definition of natural spaces. But unlike Brouwer, we are also engaged in achieving computational efficiency, as well as establishing links between \class, \intu, \russ\ and \bish\ (and formal topology).

\parr

An example of a continuous function space which cannot be represented as a natural space is the space of continuous functions from Baire space to itself (see \VelThe, copied in \WaaNat). Still there is a subset \Fun\ of Baire space \NNt\ such that every $\alpha\inn\Fun$ codes a natural morphism from Baire space to itself, and every natural morphism from Baire space to itself is coded by some $\alpha\inn\Fun$.
 
\sbsc{Basic-open spaces and basic neighborhood spaces}\label{basicneighbor}\hspace*{-2.5pt}
Basic dots do not always represent an open set, or even a neighborhood in the apartness topology.\footnote{In contrast to formal topology, where one only works with opens.}\ Still, so-called `basic-neighborhood spaces' are fundamental, especially in the context of metric spaces. In \class, \intu\ and \russ\ every continuous function from a natural space to a basic neighborhood space \Vnatt\ can be represented by a natural morphism. The idea is to look at basic dots $a$ which are neighborhoods, meaning $\closr{a}$ contains an inhabited open $U$.

\defi
Let \Vnatt\ be a natural space, with corresponding pre-natural space \Vprenatt. Let $a$ be a basic dot, and let $x\inn\hattr{a}$. Then $a$ is called a \deff{basic (open) neighborhood}\ of $x$ iff $\closr{a}$ is a neighborhood of $x$ (resp. $\closr{a}$ is itself open). Now \Vnatt\ is called a \deff{basic-open space}\ iff $\closr{a}$ is open for every $a\inn V$. \Vnatt\ is called a \deff{basic neighborhood space}\ iff \Vnatt\ is isomorphic to a basic-open space. 
\edefi

\rem
`Basic-open space' is not a `natural' property, meaning that it is not necessarily preserved under isomorphisms (see \ref{defrealnat}, where \Rprenatot\ and $(\Rrat,\aprtR,\leqcR)$ are isomorphic, yet only \Rprenatot\ is basic-open). So we `naturalize' the concept `basic-open space' to `basic neighborhood space', which then trivially is preserved under isomorphisms. 
\erem

If \Vnatt\ is a basic neighborhood space derived from \Vprenat, then $V$ contains a neighborhood basis for the natural topology. The converse does not hold in \class: see \WaaNat, A.2.5.
\parr
The prime example of a basic neighborhood space is a basic-open space where the basic dots represent open sets. We put forward the main theorem, that in \class, \intu\ and \russ\ continuous functions from a natural space to a basic neighborhood space can be represented by a natural morphism. We only need a Lindel\"of property, which follows from \BDD\ (Bar Decidable Descent, `Every bar descends from a decidable bar') which is in the common core of \class, \intu\ and \russ\ (see \WaaArt). One could see \NTop\ and its morphisms as a way of incorporating \BDD\ in the definitions, to make it accessible for \bish. Later we show that every complete metric space has a basic-open representation (and therefore in \class, \intu\ and \russ\ by the corollary below a unique representation (up to isomorphism) as a basic neighborhood space, see \ref{sepmetnat}).

\thm (in \class, \intu\ and \russ; using \BDD)

Let $f$ be a continuous function from a natural space \Vnatot\ to a basic neighborhood space \Wnattot.  Then there is a natural morphism $g$ from \Vnatot\ to \Wnattot\ such that for all $x$ in \Vcalt: $f(x)\equivvto g(x)$.

\ethm 

\prf
The not so easy proof is given in \WaaNat, A.3.1.
\eprf

\crl (in \class, \intu\ and \russ)
If \Vnatot\ and \Wnattot\ are two homeomorphic basic neighborhood spaces, then they are isomorphic.
\ecrl

\rem
The theorem suggests that from a \bish\ point of view, the concept of `natural morphism' adequately captures the notion of continuous function (under the usual topological definition). To capture the metric property `uniformly continuous on compact subspaces' we can define `inductive morphisms'. The required `genetic' induction is theoretically advantageous but unnecessary in computational practice since all `reasonably' occurring morphisms will be inductive. See also our final discussion in paragraph \ref{repfoun}.
\erem

\sbsc{Complete separable metric spaces are natural}\label{sepmetnat}
Every complete separable metric space is homeomorphic to a natural space. Therefore every separable metric space is homeomorphic to a subspace of a natural space. Some key examples of spaces which can be constructed as a natural space are \Nt, \Rt, \Ct, the complex p-adic numbers \Cpt, \RNt, Baire space, Cantor space, Hilbert space \Hit, and every Banach space. We prove slightly more, because of our interest in different representations of complete metric spaces: 

\thm
Every complete separable metric space \xdt\ is homeomorphic to a basic-open space \Vnatt. 
\ethm

\prf
The rough idea is simple: for a separable metric space \xdt\ with dense subset \anninnt, let for each $n,\sinn$ a basic dot be the open sphere $B(a_n,2^{-s})=\{x\inn X\midd d(x, a_n)\smlr 2^{-s}\}$.  Then we have an enumerable set of dots $V$ by taking $V=\{B(a_n, 2^{-s})\midd n,\sinn\}$. The only trouble now is to define \aprt\ and \leqct\ constructively, since in general for $n, m$ and $s, t$ the containment relation $B(a_n, 2^{-s})\subseteqq B(a_m,2^{-t})$ is not decidable.  We leave this technical trouble, which can be resolved using \aczo\ (countable choice), to \WaaNat, A.3.2. (The strategy is also reproduced in \ref{prfmetcompfan}). 
\eprf

\crl 
In \class, \intu\ and \russ\ the following holds:
\be
\item[(i)]
A continuous function $f$ from a natural space \Wnattot\ to a complete metric space \xdt\ can be represented by a morphism from \Wnattot\ to a basic neighborhood space \Vnatt\ homeomorphic to \xdt, by theorem \ref{basicneighbor}. 
\item[(ii)] A representation of a complete metric space as a basic neighborhood space is unique up to isomorphism.
\ee

In \bish\ the following holds:
\be
\item[(iii)] If \xdt\ and \Vnatt\ are as above in the theorem, then we can define a metric $d'$ on \Vnatt\ (see \WaaNat, 4.0.2) by defining $d'(x, y)\iz d(h(x), h(y))$ for $x,y\inn\Vcal$ and $h$ a homeomorphism from \Vnatt\ to \xdt. This metric can be obtained as a morphism from $(\Vcal\timez\Vcal, \Topaprtprod)$ (see \ref{infprod}) to \Rnatt\ by the construction of \Vnatt. We then see that the apartness topology and the metric $d'$-topology coincide, in other words \Vnatt\ is metrizable.
We conclude: on a well-chosen basic-neighborhood natural representation of a complete metric space, the metric topology coincides with the apartness topology.
\ee
\ecrl

\rem
\be
\item[(i)]
The construction in the proof sketch above merits a closer look, since we do not simply choose each `rational sphere' $B(a_n, q), q\inn\Q$ to be a basic dot. Yet for \Rt\ and its corresponding natural space \Rnatt, choosing all closed rational intervals works fine. We cannot guarantee in the general case \xdt\ however, that by taking $V=\{B(a_n, q)\midd q\inn\Q,\ninn\}$ we end up with a natural space \Vnatt\ which is homeomorphic to \xdt. We do know that for $X=\Cp$, taking $V=\{B(a_n, q)\midd q\inn\Q,\ninn\}$ gives us a \Vnatt\ which contains `more' points than \Cpt. In \WaaNat, A.2.1 we detail this nice example of a non-archimedean metric natural space.

\item[(ii)]
For most applied-computational purposes, a basic neighborhood representation of a complete metric space seems the best option. We believe that for \Rt, the representation \sigRt\ which we define in the following sections is a good choice for computational purposes also. Our definition of \sigRt\ and \leqc-morphisms matches the recommendations in \BauKav\ for efficient exact real arithmetic.

\item[(iii)] That the metric topology coincides with the apartness topology on (a well-chosen basic-neighborhood representation of) a complete metric space, allows for theoretical and computational simplification.
\ee
\erem

\sbsc{Metrizability of natural spaces}\label{metriz1}
From intuitionistic topology, we can retrieve results on the metrizability of natural spaces. With a definition of the notion `star-finitary' which closely resembles the notion `strongly paracompact', we obtain the constructive theorem that every star-finitary natural space is metrizable. 

\parr

Also, we can easily define natural spaces which are non-metrizable. Comparable to ideas from Urysohn (\Urya), in intuitionistic topology one finds spaces with separation properties `$T_1$ but not $T_2$' and `$T_2$ but not $T_3$' (see \WaaThe).  These spaces can be transposed directly to our setting. 

\parr

However, a different class of non-metrizable natural spaces arises when we look at direct limits in infinite-dimensional topology. As an example, in \WaaNat\ we show that the space of `eventually vanishing real sequences' (which is the direct limit of the Euclidean spaces $(\R^n)_{n\in\N}$) can be formed as a non-metrizable natural space.

\sbsc{(In)finite products are natural}\label{infprod}
The basic idea to arrive at the natural product of \texttt{a)} a finite sequence \texttt{b)} an infinite sequence of natural spaces is simple. Just take the Cartesian product of the basic dots involved, and define an appropriate \leqc-relation and \aprt-relation for this product. This is detailed in \WaaNat, the technical details are however less easy than one might expect. This is due again to the need to be theoretically complete. In practice the basic idea works fine, since it suffices for all perfect spaces.  

\sectionn{Applied math intermezzo:\hspace*{-1.5pt} Hawk-Eye, binary and decimal reals}\label{appmath}

\sbsc{Hawk-Eye}\label{hawkeye}
We can now discuss an interesting application of mathematics, in the world of professional tennis. In 2006 the multicamera-fed decision-support system \hwki\ was first officially used to give players an opportunity to correct erroneous in/out calls. \hwki\ uses ball-trajectory data from several precision cameras to calculate whether a given ball was \textsc{in}: `inside the line or touching the line' or \textsc{out}: `outside the line'. \hwki\ is now widely accepted, for decisions which can value at over \$100,000. 

\parr

The measurements of the cameras can be seen as the `dots' or `specks' that we used for illustration in our introduction. Software of \hwki\ must in some way run on these dots. The interesting thing is that \hwki\ does not have the feature of a \textsc{let}: `perhaps the ball was in, perhaps the ball was out, so replay the point'. From this and our work so far we derive: 

\clmm
\hwki, irrespective of the precision of the cameras, will systematically call \textsc{out}\ certain balls which are measurably \textsc{in}\ or vice versa.
\eclmm

\parr

The claim is not per se important for tennis. \hwki\ admits to an inaccuracy of 2-3 mm, and under this carpet the above claim can be conveniently swept (still, one sees `sure' decisions where the margin is smaller). \hwki's inaccuracy is usually blamed on inaccuracy of the camera system. But regardless of camera precision we cannot expect to solve the topological problem that there is no natural morphism from the real numbers to a two-point natural space $\{\mbox{\textsc{in}}, \mbox{\textsc{out}}\}$ which takes both values \textsc{in}\ and \textsc{out}. And our recommendation to \hwki\ is to introduce a \textsc{let} feature. Combined with a finer apartness this allows for an elegant solution, see \WaaNat, A.2.0 for a more detailed description.

\sbsc{Binary, ternary and decimal real numbers}\label{binarydecimal}
We next consider morphisms from \Rnat\ to the (natural) binary real numbers \Rbint\ and decimal real numbers \Rdect. These morphisms reveal the topology behind different representations of the real numbers on a computer, and transitions between these representations. For simplicity we discuss mainly \Rbint, since the situation with \Rdect\ is completely similar. For some purposes also the ternary real numbers \Rtert\ come in handy. 

\defi
We first put $\Rratbin\!\isdef \{\maxdotR\}\cupp\{[\frac{n}{2^{m}}, \frac{n+1}{2^{m}}]\midd n\inn\Z, \minn\}$. Similarly, let $\Rratter\!\isdef \{\maxdotR\}\cupp\{[\frac{n}{3^{m}}, \frac{n+1}{3^{m}}]\!\mid n\inn\Z,\minn\}$ and $\Rratdec\isdef \{\maxdotR\}\cupp[\frac{n}{10^{m}}, \frac{n+1}{10^{m}}]\midd n\inn\Z, \minn\}$. 

Then 
$\Rbin=(\Rratbin, \aprtR, \leqcR)$ is the natural space of the \deff{binary real numbers}. Similarly we form the corresponding natural spaces \Rtert\ and \Rdect\ of the \deff{ternary}\ and
\deff{decimal real numbers}.

\parr
Put $\zorratbin\isdef \{[\frac{n}{2^{m}}, \frac{n+1}{2^{m}}]\midd n,\minn \midd n\smlr 2^{m}\}$, $\zorratter\isdef \{[\frac{n}{3^{m}}, \frac{n+1}{3^{m}}]\mid n,\minn \mid n < 3^{m}\}$ and $\zorratdec\isdef \{[\frac{n}{10^{m}}, \frac{n+1}{10^{m}}]\midd n,\minn \midd n\smlr 10^{m}\}$ to form the corresponding natural spaces \zorbint, \zortert\ and \zordect, each with the same maximal dot $[0,1]$ denoted by $\maxdot_{\zor}$. 

\parr

As a partial order, $(\Rratbin, \leqcR)$ is a tree. The natural binary reals correspond to real numbers $x$ which can be given as a \deff{binary expansion} $x=(\mino)^s\cdott\Sigma_{\ninn} a_n\cdott 2^{-n+m}$, where $s\inn\{0, 1\}$, $m\inn\N$ and $a_n\inn\{0, 1\}$ for all \ninn, such that $m\bygr 0$ implies $a_0\notiz 0$. We call $s$ the \deff{sign}\ and write $s\iz +, -$ for $s\iz 0, 1$ respectively. We call $m$ the \deff{binary point place}. Then the $(a_n)_{\ninn}$ are the \deff{binary digits} in this binary expansion of $x$, and we write $x\iz (s)\,a_0\, a_1 \ldots a_{m}{\rm \bol{.}} a_{m+1} \ldots$. Notice the \deff{binary point}\ that we write between $a_m$ and $a_{m+1}$ to denote the binary point place.

Replacing `binary, 2' with `ternary, 3' and `decimal, 10' respectively, we obtain the similar definitions for \Rtert\ and \Rdect. 
\edefi

Classically every real number $y$ has an equivalent binary expansion, but in computational practice and in constructive mathematics this is not the case (see \mbox{e.g.}\ \GNSWcom\ for a thorough discussion). So with \Rbint, \Rtert\ and \Rdect\ we in practice obtain \emph{different} representations of the real numbers. We wish to shed some light on the natural topology involved in the (im)possible transition from one such representation to another.

\sbsc{Morphisms to and from the binary reals}\label{morfbinreals}
It turns out that a morphism $f$ from \Rnat\ to \Rbint\ which is order preserving ($x\leqqR y$ implies $f(x)\leqqR f(y)$) has to be `locally constant' around the $f$-originals of the rationals $\{\frac{k}{2^{m}}\midd k\inn\Z,\minn\}$. For these rationals the binary expansion has two alternatives (\mbox{e.g.}\ for $1$ both $0.111\ldots\equivv 0+1\cdott 2^{-1}+1\cdott 2^{-2}+1\cdott 2^{-3}+\ldots$ and $1.000\ldots\equivv 1+0\cdott 2^{-1}+0\cdott 2^{-2}+0\cdott 2^{-3}+\ldots$ are binary representations). Since these binary rational numbers lie dense in \Rt, there can be no injective morphism from \Rnatt\ to \Rbint\ (notice that any injective morphism $f$ from \Rnatt\ to \Rnatt\ is either order preserving, or order reversing in which case a similar argument for local constancy obtains). But this does not mean that all morphisms from \Rnatt\ to \Rbint\ are constant.

\parr

The well-known \deff{Cantor function}\ \fcant\ (also known as `the devil's staircase') is an example of a non-constant natural morphism from \zort\ to \zorbint. The Cantor function is most easily described as a refinement morphism from \zortert\ to \zorbint, but also can be given as a trail morphism on \zort, see \WaaNat, A.2.2.

\parr

We now have an example in \class\ of a continuous function between natural spaces which cannot be represented by a morphism. In \class, the identity is a homeomorphism from \Rnatt\ to \Rbint\ (remember that in \class\ we work with the equivalence classes, and that every real number has an equivalent binary representation). But this identity cannot be represented by a natural morphism, as we pointed out above. In the light of theorem \ref{basicneighbor}, the `reason' for this is that \Rbint\ is not a basic neighborhood space, which we can easily verify by looking at the real number $\half$. In fact, in \Rbint, of the basic dots only the maximal dot is a neighborhood of $\half$.

\rem
One can show with little effort that for $n, \minn$ the $n$-ary and $m$-ary reals are \leqc-isomorphic. However, we believe the $n$-ary reals can only be identically embedded in the $m$-ary reals if there is a $b\geqq 1$ in \Nt\ such that $m$ divides $n^b$ (for an identical embedding $f$ we have $f(x)\equivvR x$ for all $x$). This gives a natural-topological classification of the different $n$-ary real numbers. Furthermore, the $n$-ary reals are not constructively closed under basic arithmetic operations. Decimal representation then would seem a poor computational choice (but see \ref{disclaimer}).
\erem

\sectionn{Natural Baire space and natural Cantor space}\label{introBaire0}

\sbsc{Introduction to natural Baire space}\label{introBaire2}
We will show that natural Baire space $\Bnat=\NNnat$ is a universal natural space, meaning that every natural space \Vnatt\ is the image of natural Baire space under some natural morphism from \Bnatt\ to \Vnatt.

\parr

Natural spaces therefore correspond to quotient topologies of Baire space which are derived from a $\Sigma^0_1$-apartness (to see that this class is larger than the class of Polish spaces, it suffices to see that some of these quotient spaces are non-metrizable). The obvious representation of Baire space is well-suited for computational purposes. In the remainder of the paper we study how to transfer this property to (suitable representations of) other natural spaces.

\parr

Quite some work has already been done in intuitionistic topology, with Baire space as fundament. In \WaaNat\ it is shown that the setting of natural spaces mirrors Brouwer's setting in many ways. This means that we can simply transpose many intuitionistic results. But we take a neutral constructive approach, and do not use any specific classical or intuitionistic axioms. However, in developing the theory we freely use the axioms of countable choice \aczo\ and dependent choice \dco, which are generally accepted as constructive.

\sbsc{Natural Baire space}\label{natBaire} 
The definition of Baire space as natural space is simple. Its set of basic dots is \Nstart, the set of all finite sequences of natural numbers (representing the basic clopen sets of Baire space). We use the definition also to relate natural Baire space to usual Baire space.

\defi 
Let \Nstart\ be the set of all finite sequences of natural numbers. For $a=a_0, \ldots, a_i$, $b\iz\ b_0, \ldots, b_j\inn\Nstar$ the concatenation $a_0,\ldots, a_i, b_0\ldots, b_j$ is denoted by $a\star b$. Define: $b\leqcom a$ iff there is $c$ such that $b=a\star c$. Define: $a\aprtom b$ iff $a\nleqcom b$ and $b\nleqcom a$. 

\parr

Then \Bprenatt\ is a pre-natural space, with the empty sequence as maximal dot, which we also denote \maxdotBt\ or simply \maxdott. Its corresponding natural space \Bnattopt\ we call \deff{natural Baire space}. We also write \NNnatt\ for \Bnatt.

\parr

Given $\alpha\inn\NN$ and \minnt, we write $\alfstr(m)$ for the finite sequence $\alpha(0), \ldots, \alpha(m-1)$ consisting of the first $m$ values of $\alpha$. Notice that $\alfstr(m)$ is an element of \Nstart, so the sequence $\alfstr\iz \alfstr(0), \alfstr(1), \ldots $ is a point in \Bnat.

Conversely, for a point $p\inn\Bnat$, there is a unique sequence $\alpha\inn\NN$ such that $p\equivvom\alfstr$. We write $p^*$ for this unique $\alpha$, giving that $p\equivvom\overline{p^*}$ for $p\inn\Bnat$ and $\alpha\iz\alfstr^*$ for $\alpha\inn\NN$.
\edefi

\thm
\Bnattop\ is homeomorphic with \Bairut.
\ethm

\prf
We leave it to the reader to verify that the function $\alpha\rightarrow\alfstr$ from $\NN$ to \Bnatt\ defined above is a homeomorphism, with inverse $p\rightarrow p^*$ (also defined above).
\eprf 

\sbsc{Natural Cantor space}\label{defcantor}
We first define the notion `natural subspace', since in natural Cantor space we have a prime example.

\defi Let \Vnatt\ be a natural space derived from \Vprenatt. Let \Wt\ be a countable subset of \Vt, then \Wprenatt\ is a pre-natural space, with corresponding set of points \Wcalt. If \Wnatt\ is a natural space (see def.\,\ref{defnatspace}), then we call \Wnatt\ a \deff{natural subspace}\ of \Vnatt\ iff in addition \Wnatt\ as a natural space coincides with \Wnatt\ as a topological subspace of \Vnatt\ (in the subspace topology `$U\subseteq\Wcal$ is open' is defined thus: there is an open $U'\subseteq\Vcal$ such that $U\iz U'\capp\Wcal$).

\parr

Let \zostart\ be the set of finite sequences of elements of \zot. Now \deff{natural Cantor space}\ is the natural subspace \Cnattopt\ of natural Baire space formed by the pre-natural space \Cprenatt\ and its set of points \Cnatt.
\edefi

\rem Natural Cantor space is homeomorphic to usual Cantor space, and corresponds directly to Brouwer's fan \sigtot. 
\erem

From now on, when the context is clear we will simply say `Baire space' and `Cantor space' and omit the extra word `natural'.

\sectionn{Lattices, trees and spreads}\label{lattrespre}

\sbsc{Lattices and posets of basic dots}\label{latopen}
In topology, the open sets form a lattice structure under the inclusion relation. This structure is often exploited in various ways, for instance in domain theory and formal topology. Since we have the added apartness, we can disregard meet and join operations and focus simply on the partial-order properties (of the `poset' of opens). Our basic dots in general need not form a lattice, but their partial-order properties play an important role. We now go into these partial-order properties in more detail.

\sbsc{Trees and treas}\label{treas}
For Baire space the poset of basic dots \Nleqcomt\ forms a countable tree. That is: for any dot $a=a_0,\ldots, a_{n-1}\inn\Nstar$, there is a \emph{unique} finite trail of immediate successors/predecessors from \maxdotBt\ to $a$. (Therefore any \precc-trail between dots is finite, and also the successor/predecessor relationship is decidable.). We cannot achieve this elegance for any natural space, but we can show that any natural space \Vnatot\ is isomorphic to a natural space \Wnattot\ where \Wleqctot\ equals \Nleqcomt. Or more practical: where \Wleqctot\ is a full subtree of \Nleqcomt, definition follows.

\parr

This means that we could limit ourselves to natural spaces \Vnatt\ where \Vleqct\ is (a full subtree of) \Nleqcomt. But we see two possible downsides to limiting ourself to  \Nleqcomt. One downside is that for many natural spaces, we would have to replace our original basic dots with elements of \Nstart, which can be a tedious encoding.\footnote{Basic dots always derive from \Nt, but still we prefer to write `$[\half,1\half]$' etc.}\ The other downside we see is that we often have to introduce duplicate copies of basic dots, in order to obtain a tree. These duplicates multiply in number with every refinement step, which seems hardly efficient when building actual implementations.

\parr
Therefore we propose the compromise notion of a `trea'. One can think of a trea as being a tree wherein certain of the branches have been neatly glued together in a number of places. An important example is the set of the lean dyadic intervals $\sigR\iz(\{\maxdotR\}\cupp\{[\frac{n}{2^{m}}, \frac{n+2}{2^{m}}]\midd n\inn\Z, \minn\}, \leqcR)$. A more precise characterization of a trea: a countable \precc-directed acyclic graph with a maximal element, where for each node there are finitely many immediate-predecessor trails to the maximal element, all of the same length. Another characterization: a countable \leqct-poset with a maximal element where each point has finitely many immediate-predecessor trails to the maximal element, all of the same length.

\defi   
Let \Vnatt\ be a natural space, with corresponding \Vprenat,  and let \Wnatt\ with corresponding \Wprenatt\ be a natural subspace of \Vnatt\ (so $W\subseteq V$). Let $a\precc c$ in $V$.  
\be
\item[(i)]
We say that $a$ is a \deff{successor} of $c$ in \Vleqct\ (notation $a\suczV c$, or simply $a\sucz c$ if the context is clear) iff for all $b\inn V$, if $a\precc b\leqc c$ then $b\iz c$. A sequence  $b_0\predcz\ldots\predcz b_n$ in $V$ is called a
\deff{\sucz-trail of length $n$ from $b_0$ to $b_n$ in \Vleqct}. For $b\inn V$ we put $\suczV(b)\isdef\{d\inn V\midd d\sucz b\}$, and simply write $\suczz(b)$ when the context is clear.
\item[(ii)] 
\Vleqct\ is called a \deff{trea} iff for every $a\inn V$ the set $\{b\inn V\midd a\leqc b\}$ of predecessors of $a$ is finite (then the successor relation \suczt\ is decidable, and for every $a\inn V$ there are finitely many \sucz-trails from \maxdott\ to $a$) and in addition there is an integer $\grdd(a)\inn\N$ such that every \sucz-trail from \maxdott\ to $a$ has length $\grdd(a)$.
\item[(iii)]
Now let \Vleqct\ be arbitrary, where \Wleqct\ is a tree (trea), then we say that \Wleqct\ is a \deff{subtree (subtrea)} of \Vleqct.
\item[(iv)]
Let \Vleqct\ be a tree (trea), and \Wleqct\ a subtree (subtrea). We then call \Wleqct\ a \deff{full subtree (subtrea)} of \Vleqct\ iff $b\suczW d$ implies $b\suczV d$ for all $b, d\inn W$. (Then each \suczW-trail in \Wleqct\ is a \suczV-trail in \Vleqct).

\ee
\edefi

Treas behave just like trees (and any tree is a trea). Most of the spaces of interest that we mentioned so far (see \ref{natlarge}) have an intuitive representation as a natural space \Vnatt\ where \Vleqct\ is a trea.

\xam
For the natural real numbers \Rnatt\ we can easily indicate an isomorphic subspace \signatR\ with corresponding pre-natural space \sigprenatRt, where \sigleqcRt\ is a trea: 

\parr

$\sigprenatR\isdef(\{\maxdotR\}\cupp\{[\frac{n}{2^{m}}, \frac{n+2}{2^{m}}]\midd n\inn\Z, \minn\}, \aprtR, \leqcR)$.

\parr

Our examples in \ref{binarydecimal}\ should show why we cannot hope to find an isomorphic subspace $(\Vcal,\TopRaprt)$ where \Vleqct\ is a tree (!).
\exam

\sbsc{Spreads and spraids}\label{spraids}
The previous example illuminates a bridge towards intuitionistic terminology, which we give in the following definition:

\defi
Let \Vnatt\ be a natural space, with corresponding \Vprenat,  and let \Wnatt\ with corresponding \Wprenatt\ be a \deff{decidable} natural subspace of \Vnatt\ (meaning $W$ is a decidable subset of $V$). 
\be
\item[(i)]
We call \Vnatt\ a \deff{spread (spraid)} iff \Vleqct\ is a tree (trea) and each infinite \precc-trail defines a point. Then we call \Wnatt\ a \deff{subspread (subspraid)} of \Vnatt\ iff \Wleqct\ is a full subtree (subtrea) of \Vleqct.
\item[(ii)]
We call \Vnatt\ a \deff{Baire spread} iff \Vnatt\ is a subspread of Baire space.  
\ee
 
By extension, \Vnatt\ is \deff{spreadlike}\ iff there is an isomorphism between \Vnatt\ and a spread. 
\edefi

\xam
Important basic examples of subspraids are obtained as follows. For \Vnatt\ a spraid and $a$ in $V$, one easily sees that $\Vsuba\iz\{b\inn V\midd b\leqc a\}\iz\{a\}_{\!\leqc}$ determines a subspraid of $V$ if we put its maximal dot as $\maxdot_a\iz a$.
\exam

\sectionn{Universal spaces and efficient representations}\label{uninat}

\sbsc{Baire space is universal}\label{Baireuni}
Baire space is a universal natural space, by which we mean that each natural space can be seen as the image of Baire space under a natural morphism. In other words: every natural space is spreadlike. In this article we look to use this result for computational efficiency. Theoretically, on a meta-level it gives us a direct correspondence with many important intuitionistic results. See \WaaThe, in which apartness topology is developed in \intu. 

\thm
Every natural space is spreadlike. In fact, every natural space \Vnatt\ is isomorphic to a spread \Wnatt\ whose tree is $(\Nstar, \leqcom)$.
\ethm

\crl
\be
\item[(i)] Let \Vnatt\ be a natural space, then there is a surjective \leqc-morphism from Baire space to \Vnatt. (`Baire space is a universal spread', `every natural space is the natural image of Baire space', `every natural space is a quotient topology of Baire space'). 
\item[(ii)] If \Vnatt\ is a basic-open space (see definition \ref{basicneighbor}) then \Vnatt\ is isomorphic to a basic-open spread \Wnatt\ whose tree is $(\Nstar, \leqcom)$. 
\ee
\ecrl

\prf
See \WaaNat, A.3.4.
\eprf

The corollary gives the equivalent picture that each natural space \Vnatt\ with corresponding pre-natural space \Vprenatt\ is in fact nothing but a pre-apartness \aprtVt\ on \Nstart\ which respects \aprtomt\ and \leqcomt. To define \aprtVt\ we only have to `pull back' the decidable relation \aprt\ using the given surjective morphism $f$ thus: for $a, b\inn\Nstar$ put $a\aprtV b$ iff $f(a)\aprt f(b)$ (then $a\aprtV b$ implies $a\aprtom b$).

\parr

An ideal situation which avoids encoding arises whenever a natural space \Vnatt\ contains a subspraid on which the identity is an isomorphism with the whole space. Then from the often vast partial-order universe of \Vleqct\ we can restrict ourselves to a subtrea. We give the important example of the real numbers below, where the isomorphic subspace is a spraid. We believe this to be the most common setting for natural spaces. In the uncommon case that we cannot find an isomorphic subspace which is a spraid, we can always find an isomorphic spread.  

\xam
Looking at the natural real numbers \Rnatt, we can easily indicate an isomorphic subspace which is a spraid as in example \ref{treas}. Put 

\parr

$\sigR\isdef \{\maxdotR\}\cupp\{[\frac{n}{2^{m}}, \frac{n+2}{2^{m}}]\midd n\inn\Z, \minn\}$.

\parr

Then \sigprenatRt\ is a spraid which is an isomorphic subspace of \Rnatt. Similarly we define:

\parr

$\sigzor\isdef \{[\frac{n}{2^{m}}, \frac{n+2}{2^{m}}]\midd n, \minn\midd n\pluz 2\leqq 2^{m}, m\geqq 1\}$,

\parr
so that taking $\maxdotzor\iz [0,1]$ we get a subfann (see next paragraph)  \sigzort\ of \sigRt\ which is isomorphic to \zornatt. 
\exam

\parr
Another more involved example of a spraid arises when building the natural space $C^{\rm unif}(\zor, \R)_{\rm nat}$ of uniformly continuous functions from \zort\ to \Rt. This is sketched in \WaaNat, A.2.4, referring for details to earlier work of Brouwer.

\sbsc{Cantor space is a universal fan}\label{Cantuni}
Where Baire space is a universal spread, Cantor space is a universal fan, by which we mean that each `finitely branching' spraid can be seen as the image of Cantor space under a natural morphism:

\defi
Let \Vnatt\ be a spread (spraid) derived from \Vprenat. We call the tree (trea) \Vleqct\ \deff{finitely branching} iff for all $c\inn V$ the set $\suczz(c)\iz\{a\inn V\midd a\sucz c\}$ is finite. We call \Vnatt\ a \deff{fan (fann)} iff \Vleqct\ is a finitely branching tree (trea). By extension, \Vnatt\ is \deff{fanlike}\ iff \Vnatt\ is isomorphic to a fan.
\edefi

\thm
Let \Vnatt\ be a fann, then there is a surjective morphism from Cantor space to \Vnatt. (`Cantor space is a universal fan').
\ethm

\crl 
Every fann is fanlike. Every fanlike space is the natural image of Cantor space.
\ecrl

\prf
See \WaaNat, A.3.5.
\eprf

\sbsc{Every compact metric space is homeomorphic to a fan}\label{metcompfan}
If we define a separable metric space to be compact whenever it is totally bounded and complete (as is standard in \bish), then it is a well-known result that every compact metric space is the uniformly continuous image of Cantor space. The following theorem is therefore not surprising. It shows that we can represent a compact metric space with a fan (the theorem is not mentioned explicitly in \WaaNat, which is why we need to prove it in the appendix):

\thm
Every compact metric space is homeomorphic to a fan.
\ethm

\prf
For a given compact metric space, we need to find a suitable fan and prove that the apartness topology coincides with the (induced) metric topology. We solve this using the theory of \WaaNat, see the appendix \ref{prfmetcompfan}.
\eprf

\rem
To use compactness in general topology, for \bish\ it seems unavoidable to adopt (transfinite) inductive machinery. So-called genetic induction is developed in \WaaNat\ to enable the use of Heine-Borel properties of compact spaces for theoretical purposes. For computational practice genetic induction is unimportant since Brouwer's Fan Theorem (\FT) `always' holds (see also \ref{repfoun}).
\erem 

\sbsc{Representation of (locally) compact metric spaces}\label{reploccomp}
By the previous theorem every compact metric space can be represented by a fan. But just like the situation with the real numbers, it is often more convenient to represent a compact metric space by a fann. In computational practice it will not be difficult to find good representations. The same holds for locally compact metric spaces. They can be represented by a spraid which consists of a countable number of fanns.

\sbsc{Representation of complete metric spaces}\label{repcompmet}
For complete metric spaces which are not locally compact, we can use the property `strong paracompactness'. In \class\ a complete metric space can be represented by a `star-finite' spread. Conversely, already in \bish\ we can show that `star-finitary' natural spaces are metrizable. `Star-finite' is a generalization of `locally finitely branching' (see \WaaNat, 4.0.7--4.0.10) where each basic dot only has finitely many touching neighbors of the same length. Once again, we think that in computational practice it will not be difficult to find good star-finite representations of a given complete metric space.

\sectionn{Efficient computation on spreads and spraids}\label{spreadmorph}

\sbsc{Refinement versus trail morphisms 1}\label{refvstrail1}
We return briefly to our discussion of refinement morphisms versus trail morphisms. With spreads (which derive from a tree) there is no need for trail morphisms. In fact a spread \Vnatt\ is \leqc-isomorphic to its trail space (denoted \Vpathnatt, see \ref{deftrailmorph}).
Since Baire space is universal (\ref{Baireuni}), we could develop a fruitful theory using only spreads and refinement morphisms (as is done in \intu). 

For computational purposes, one sees that a spread representation of the real numbers is cumbersome when compared to the spraid \sigRt.
Studying refinement morphisms also on spraids therefore seems a fruitful endeavour. But when working with spraids, we in theory sometimes need trail morphisms as well. Fortunately, we can show that continuous real functions can always be represented by a refinement morphism from \sigRt\ to \sigRt. More generally: any continuous function from a spraid \Vnatt\ to the reals can be represented by a refinement morphism from \Vnatt\ to \sigRt.

\parr

We see this by first looking at proposition \ref{refvstrail2} below, which states that a trail morphism from a spraid \Vnatt\ to \sigRt\ can already be represented by a refinement morphism from \Vnatt\ to \sigRt. We then combine this with theorem \ref{basicneighbor} that every continuous function from a natural space to \sigRt\ can be represented by a morphism.

\parr

We take some time to argue that for many important spraids resembling \sigRt, a trail morphism can already be directly represented by a refinement morphism. This is especially relevant for the computational perspective, we believe. To keep our narrative comprehensive, we use representation results from \WaaNat, and refer also to \WaaNat\ for the technical details.

\sbsc{Unglueing of spraids}\label{unglspraids}
Spraids correspond to treas, which can be seen as trees in which branches are glued together in a certain regular way (see \ref{treas}). To simplify the theory, we note that any spraid \Vnatt\ can be unglued in a simple manner to a spread \Vunglnatt. We specify this in the appendix \ref{unglueing}. Unglueing a spraid \Vnatt\ amounts to adding, for each $a\inn V$, a finite number of copies of $a$ such that each \sucz-trail from \maxdott\ to $a$ is represented by one of the copies. These copies all have $\grdd(a)$ as length in $(\Vungl, \leqcstar)$.

For spraids, working with \Vunglnatt\ is more elegant than working with the trail space \Vpathnatt. If we start with a spread \Vnatt, then there is a trivial bijection between $V$ and \Vunglt, showing that spreads are already unglued.

\sbsc{Refinement versus trail morphisms 2}\label{refvstrail2}
Now we can show that for many important spraids a trail morphism can already be directly represented by a refinement morphism. We illustrate this first with \sigRt, our preferred representation of \Rt. Therefore \bisco\ real functions (uniformly continuous on closed intervals) can always be represented by a \leqc-morphism sending lean dyadic intervals to lean dyadic intervals. This paragraph corresponds (we believe) to propositions 4.2 and 8.2 in \BauKav, but we do not need Markov's Principle. The difference seems cosmetic though, and for practice even non-existent. We use a similar lazy convergence, but avoid the axiom.   See our discussion in \ref{repfoun}. 

\prp (expanded from \WaaNat)
Let \Vnatt\ be a spraid. Let $f$ be a \pthh-morphism from \Vnatt\ to \sigRt. Then there is a \leqc-morphism $g$ from \Vnatt\ to \sigRt\ such that $f(x)\equivvR g(x)$ for all $x\inn\Vcal$.
\eprp

\prf
We see $f$ as a \leqc-morphism from \Vunglnatt\ to \sigRt. 
For $c\inn\sigR$ of the form $[\frac{4s+i}{2^{t+2}}, \frac{4s+i+2}{2^{t+2}}]$ where $1\leqq i\leqq 4$ and $t\inn\N$, put $\widehat{c}\iz [\frac{s}{2^{t}}, \frac{s+2}{2^{t}}]$. For all other
$c\inn\sigR$ let $\widehat{c}\iz\maxdotR$. 
Now for $a\inn V$ there are finitely many \sucz-trails from \maxdotVt\ to $a$, say $b_0, \ldots, b_n$ where each $b_i$ is in \sigRunglt. Since the $f(b_i)$'s all touch, $\bigcapp_i\,\widehat{f(b_i)}$ is in \sigRt. We 
put $g(a)\isdef\bigcapp_i\,\widehat{f(b_i)}$. Then $g$ thus defined is a \leqc-morphism from \Vnatt\ to \sigRt\ such that $f(x)\equivvR g(x)$ for all $x\inn\Vcal$.
\eprf

\thm
Let $f$ be a \bisco\ real function (that is: uniformly continuous on closed intervals; w.l.o.g. from \Rnatt\ to \Rnatt\ by thm.\,\ref{defrealnat}). Then there is a \leqc-morphism  $g$ from \sigRt\ to \sigRt\ such that $f(x)\equivvR g(x)$ for all $x\inn\R$.
\ethm

\prf
In \WaaNat, 3.3.3, it is proved that every \bisco\ real function can be represented by an (inductive) morphism (and vice versa). Now apply the previous proposition. The next corollary is a not-too-difficult generalization.
\eprf

\crl
Let $f$ be a \bisco\ function from $\R^n$ to $\R^m$. Then there is a \leqc-morphism  $g$ between the corresponding lean spraid representations (using lean dyadic $n, m$-dimensional boxes) such that $f(x)\equivv g(x)$ for all $x\inn\R^n$.
\ecrl

\sbsc{Efficient representation of continuous functions}\label{repcontreal}
This paragraph is the finale of our computational narrative. We start with continuous real-valued functions. By the previous proposition we obtain:

\thm (in \class, \intu, \russ; using \BDD) 
Let \Vnatt\ be a spraid, and let $f$ be a continuous function from \Vnatt\ to \Rnatt. Then there is a \leqc-morphism  $g$ from \Vnatt\ to \sigRt\ such that $f(x)\equivvR g(x)$ for all $x\inn\Vcal$.
\ethm

\prf Just combine proposition \ref{refvstrail2} with theorem \ref{basicneighbor}.
\eprf 

\crl (in \class, \intu, \russ) In particular, every continuous function from \Rt\ to \Rt\ can be represented by a \leqc-morphism from \sigRt\ to \sigRt.
\ecrl

\parr

The theorem is only partly mentioned in \WaaNat, although the main ingredients are all present. In our eyes it shows a way to compute efficiently with topological spaces, within a robust theoretical framework.
For many complete metric spaces, a similar theorem to the above theorem for \sigRt\ holds. Paragraphs 3.4.3 and 4.0.9, 4.0.10 in \WaaNat\ illustrate that for many complete metric spaces, we can find efficient spraid representations. Our final question then is this. Given such a spraid \Vnatt, what extra properties would ensure that we can always represent continuous functions to \Vnatt\ with refinement morphisms, as in the case of \sigRt?

\parr
Given such a spraid, a sufficient property is that for a finite intersection of basic dots we can find a basic dot of `small enough diameter' which contains the intersection in its interior. (Comparable to $\bigcapp_i\,\widehat{f(b_i)}$ being in \sigRt, in the proof of the above proposition). For our standard basic-open complete metric spraids this property holds, but these spraids are themselves not an efficient representation.
\parr
We think that for (locally) compact spaces we can find suitable fanns. In the non-locally-compact situation, we think suitable star-finite representations can be found. See also paragraph \ref{repcompmet} and \WaaNat, 4.0.10.

\sbsc{Final discussion: representation and foundations}\label{repfoun}
Note that lean dyadic intervals can be stored efficiently, using just two numbers $n\inn\Z, \minn$ to denote $[\frac{n}{2^{m}}, \frac{n+2}{2^{m}}]$. Still, we deviate slightly from \BauKav\ (which uses all dyadic intervals, see also its discussion section 10). We are not knowledgeable enough to see whether this gives computational inefficiency compared to \BauKav, but we could switch to the spraid of all dyadic intervals if necessary.
\parr
In theorem \ref{refvstrail2}, \biscont\ is needed only for its Lindel\"of property (obtained with countable choice). Using \BDD\ we can then generalize to theorem \ref{repcontreal}. Our framework serves to clarify this type of axiomatic dependencies. An objective of \NTop\ in combination with \WaaArt\ is to aid the development of one robust framework for \bish\ (comparable to \class, \intu\ and \russ) for the different branches of mathematics. The current situation in constructive mathematics is more like patchwork, we feel. It is daunting to untangle the various interdependencies, and see what exactly are the underlying axioms, definitions and assumptions. In our eyes this hardly makes for an attractive theory. It also makes it difficult to check whether given representations confirm to (other) theoretical specifications. A good constructive framework should therefore be simple, in our (perhaps not so humble) opinion. We believe that Bishop-style pointwise mathematics is both attractive and sufficient.
\parr
In \WaaArt\ it was shown that the statement `\bisco\ functions are closed under composition' implies the Fan Theorem (\FT). Moreover, it was shown that this situation cannot be remedied by a simple change of definition, unless one sacrifices the `uniformly continuous on compact subspaces' condition. Unfortunately, an effect of \WaaArt\ seems to have been a steering away from the pointwise approach, in favour of pointfree topology and domain theory (or comparable). Often the flaw in \biscont\ (already mentioned in \WaaThe) is given explicitly as one reason to favour the (inductive) pointfree approach. Yet in \WaaNat, 3.4.0 we show that the pointwise situation is not entirely remedied by switching to inductive definitions. Inductive definitions tend to obfuscate that the problem lies with \russ, by excluding valid parts of \russ\ rather silently. Better to deal with it explicitly, we think, and stay in the true spirit of \bish. Else, adopting the induction axiom \BT\ (Brouwer's Thesis, which implies both \FT\ and \BDD, see \WaaArt) seems a more elegant option. 
\parr
If we adopt \BT, we can directly translate important intuitionistic results to \NTop\ and dispense with much of the inductive machinery, for elegance and simplicity. Since \BT\ is also valid in \class, this provides a simple way for a classical mathematician to appreciate intuitionistic results. Moreover, we concur with \BauKav\ that in computational practice \BT\ always holds.\footnote{Unless one is explicitly implementing some Kleene-tree based recursive counterexample, which even seems hard to do.}\ Brouwer's meta-analysis of how we can attain infinite knowledge (only through induction) looks as valid in its context as Church's Thesis. We therefore hope that the resistance to \BT\ which started with Bishop (who called Brouwer's theory of the continuum semimystical) will dwindle in the future.
\parr
The same holds for Markov's Principle (\MP). In \BauKav\ it is used, but we believe this use to be inessential. In \WaaArt\ it is argued that \MP\ is a form of induction comparable to natural induction over \Nt. These matters are worthy of attention, we think, to build a robust framework for \bish. 
\parr
Such a framework seems necessary, for more than one reason.
For instance, the status of the works \BisFou, \BisBri\ and \BriCon\ (analysis) is unclear due to the difficulty with \biscont\ mentioned above. Formal topology, while resolving this difficulty for a pointfree setting, seems unsuited for pointwise analysis. We think \NTop\ gives a way to retain pointwise analysis and restore many earlier \bish\ results.

\sectionn{Acknowledgements and bibliography}

\sbsc{Acknowledgements}
Wim Couwenberg came up with the basic idea for Natural Topology, and played an indispensable sparring role in many discussions. For this article, we also rely on the work done in \BauKavto\ and \BauKav. \NTop\ resembles other developments, notably intuitionistic topology, domain theory and formal topology. In \CoqFor\ the combination of refinement and apartness is already suggested. Further acknowledgements and 
historical comments can be found in \WaaNat. All (inevitable) omissions, oversights and mistakes are the author's. 

\sbsc{Bibliography}

\blit{\TroDal\hspace*{.5cm}}

\bitem{[AczCur2010]}{P. Aczel and G. Curi}
             {On the T1 axiom and other separation properties in constructive point-free and point-set topology}
             {Annals of Pure and Applied Logic vol. 161, iss. 4, pp. 560-569}{2010}{\AczCur}

\bitem{[BauKav2008]}{A. Bauer and I. Kavkler}
             {Implementing Real Numbers With RZ}
{Electronic Notes in Theoretical Computer Science 202, pp. 365-384}{2008}{\BauKavto}

\bitem{[BauKav2009]}{A. Bauer and I. Kavkler}
             {A constructive theory of continuous domains suitable for implementation}
{Annals of Pure and Applied Logic vol. 159, iss.1-3, pp. 251-267}{2009}{\BauKav}

\bitem{[Bee1985]}{M. Beeson}
           {Foundations of Constructive Mathematics}
           {Sprin\-ger-Verlag, Berlin Heidelberg}{1985}{\BeeFou}

\bitem{[Bis1967]}{E. Bishop}
           {Foundations of Constructive Analysis}
           {McGraw-Hill, New York}{1967}{\BisFou}

\bitem{[BisBri1985]}{E. Bishop and D.S. Bridges}
                      {Constructive Analysis}
                      {Sprin\-ger-Verlag, Berlin Heidelberg} {1985}{\BisBri}

\bitem{[Bri1979]}{D.S. Bridges}
 	{Constructive Functional Analysis}
                        {Pitman, London} {1979}{\BriCon}

\bitem{[BriRic1987]}{D.S. Bridges and F. Richman}
	{Varieties of constructive mathematics}
	{London Math. Soc. Lecture Notes no. 93, Cambridge University Press}{1987}{\BriRic}

\bitem{[BriV\^{i}\c{t}2006]}{D.S. Bridges and L.S. V\^{i}\c{t}\u{a}}
	{Techniques of constructive analysis}
	{Universitext, Springer Science+Business Media}{2006}{\BriVit}

\bitem{[BriV\^{i}\c{t}2011]}{D.S. Bridges and L.S. V\^{i}\c{t}\u{a}}
	{Apartness and Uniformity}
	{Springer-Verlag, Berlin Heidelberg}{2011}{\BriVitthr}

\bitem{[Bro1922]}{L.E.J. Brouwer}
              {Besitzt jede reelle Zahl eine Dezimalbruch-Entwickelung?}
              {Math. Annalen 83, 201-210}{1922}{\BroDez}

\bitem{[Bro1975]}{L.E.J. Brouwer}
              {Collected works}
	{(vol. I, II) North-Holland, Amsterdam}{1975}{\BroCol}

\bitem{[Coq1996]}{T. Coquand}
              {Formal topology with posets}
	{preprint available from www.cse.chalmers.se/~coquand/alt.ps}{1996}{\CoqFor}

\bitem{[ColEva2008]}{H. Collins and R. Evans}
	{You cannot be serious! Public Understanding of Technology with special reference to `Hawk-Eye'.}
	{Public Understanding of Science, vol. 17, 3}{2008}{\ColEva}

\bitem{[CSSV2003]}{T. Coquand, G. Sambin, J. Smith, S. Valentini}
             {Inductively generated formal topologies}
             {Annals of Pure and Applied Logic, vol. 124, iss. 1-3, pp. 71-106}{2003}{\CSSVind}

\bitem{[FouGra1982]}{M.P. Fourman, R.J. Grayson}
             {Formal Spaces}
             {in `The L. E. J. Brouwer Centenary Symposium', North
Holland, Amsterdam, pp. 107-122}{1982}{\FouGra}

\bitem{[Fre1937]}{H. Freudenthal}
              {Zum intuitionistischen Raumbegriff}
               {Compositio Mathematica, vol. 4, pp. 82-111}{1937}{\Freu}

\bitem{[GNSW2007]}{H. Geuvers, M. Niqui, B. Spitters and F. Wiedijk}
                {Constructive analysis, types and exact real numbers}
               {Mathematical Structures in Computer Science, vol.17, iss. 1, pp 3-36}{2007}{\GNSWcom}

\bitem{[KalWel2006]}{I. Kalantari and L. Welch}
                     {Larry Specker's theorem, cluster points, and computable quantum functions}
                    {Logic in Tehran, 134-159, Lect. Notes Log., 26, Association for Symbolic Logic}{2006}{\KalWel}

\bitem{[Kle\-Ves1965]}{S.C. Kleene and R.E. Vesley}
                     {The Foundations of Intuitionistic Mathematics --- especially in relation to recursive functions}
                    {North-Holland, Amsterdam}{1965}{\KleVes}

\bitem{[KreSpi2013]}{R. Krebbers and B. Spitters}
                     {Type classes for efficient exact real arithmetic in Coq}
                    {Logical Methods in Computer Science Vol. 9(1:01)2013, pp. 1-27}{2013}{\KreSpi}

\bitem{[KunSch2005]}{D. Kunkl and M. Schr\"{o}der}
                     {Some Examples of Non-Metrizable Spaces Allowing a Simple Type-2 Complexity Theory}
                    {Electronic Notes in Theoretical Computer Science, vol. 120, 111-123}{2005}{\KunSch}

\bitem{[Kus1985]}{B.A. Kushner}
	{Lectures on Constructive Mathematical Analysis}
	{American Mathematical Society, Providence, R.I.}{1985}{\KusLec}

\bitem{[M-L\"{o}f1970]}{P. Martin-L\"{o}f}
              {Notes on Constructive Mathematics}
             {Almqvist {\&} Wiksell, Stockholm}{1970}{\MLof}        

\bitem{[vMil1989]}{J. van Mill}
              {Infinite-dimensional Topology}
             {North-Holland, Amsterdam}{1989}{\Mill}        

\bitem{[Pal2005]}{E. Palmgren}
	{Continuity on the real line and in formal spaces}
	{in: `From Sets and Types to
Topology and Analysis' ed. L. Crosilla and P. Schuster, Oxford University Press}{2005}{\PalCon}

\bitem{[Pal2009]}{E. Palmgren}
	{From intuitionistic to formal topology: some remarks on the foundations of homotopy theory}
	{in: Logicism, Intuitionism, and Formalism - what has become of them?, Springer Netherlands, pp. 237-253}{2009}{\PalInt}

\bitem{[Sam2003]}{G. Sambin}
	{Some points in formal topology}
	{Theoretical Computer Science, vol. 305 iss. 1-3, Elsevier}{2003}{\SamFor}

\bitem{[TayBau2009]}{P. Taylor and A. Bauer}
                          {The Dedekind reals in Abstract Stone Duality}
                  {Mathematical Structures in Computer Science
vol. 19, iss. 4}{2009}{\TayBau}

\bitem{[Tro{\&}vDal1988]}{A.S. Troelstra and D. van Dalen}
                          {Constructivism in Mathematics}
                  {(vol. I, II) North-Holland, Amsterdam}{1988}{\TroDal}

\bitem{[Ury1925a]}{P. Urysohn}
                       {\"{U}ber die M\"{a}chtigkeit der Zusammenh\"{a}ngende Mengen}
                         {Mathematische Annalen vol. 94, pp. 262-295}{1925}{\Urya}

\bitem{[Ury1925b]}{P. Urysohn}
                       {Zum Metrisationsproblem}
                         {Mathematische Annalen vol. 94, pp. 309-315}{1925}{\Uryb}

\bitem{[Vel1981]}{W.H.M. Veldman}
              {Investigations in intuitionistic hierarchy theory}
               {PhD thesis, University of Nijmegen}{1981}{\VelThe}

\bitem{[Vel1985]}{W.H.M. Veldman}
              {Intu\"\i tionistische wiskunde}
                {(lecture notes in Dutch) University of Nijmegen}{1985}{\VelLec}

\bitem{[Vel2011]}{W.H.M. Veldman}
              {Brouwer's Fan Theorem as an axiom and as a
contrast to Kleene's Alternative}
                {Research report, on arxiv.org, arXiv:1106.2738v1}{2011}{\VelFan}

\bitem{[Waa1996]}{F.A. Waaldijk}
                    {modern intuitionistic topology}
                    {PhD thesis, University of Nijmegen}{1996}{\WaaThe}

\bitem{[Waa2005]}{F.A. Waaldijk}
                    {On the foundations of constructive mathematics --- especially in relation to the                      theory of continuous functions}
                    {Foundations of Science, vol. 3, iss. 10, pp. 249-324}{2005}{\WaaArt}

\bitem{[Waa2012]}{F.A. Waaldijk}
                    {Natural Topology}
                    {Research monograph (2nd ed., 1st ed.\, 2011), on arxiv.org, arXiv:1210.6288v1}{2012}{\WaaNat}

\bitem{[Wei2000]}{K. Weihrauch}
	{Computable analysis}
	{Springer Verlag Berlin}{2000}{\WeiCom}

\elit

\sectionn{Appendix: technical definitions and proofs}

\sbsc{Trail spaces and trail morphisms}\label{deftrailmorph}
Actually, a trail morphism from a natural space \Vnatt\ to another space \Wnattot\ is given by a refinement morphism from the `trail space' associated with \Vnatt, to \Wnattot. To define this trail space, we form new basic dots from finite sequences of `old' basic dots.

\defi
Let \Vnatt\ be a natural space derived from \Vprenatt. Let \ninnt, and let $a\iz a_0\geqc\ldots\geqc a_{n-1}$ be a shrinking sequence of basic dots in $V$. The \deff{\precc-trail of $a$}, notation $\asstr$, is the longest subsequence $a_0\succ\ldots\succ a_s$ of $a$.

\parr

For $p\iz p_0, p_1,\ldots$ in \Vcalt\ and \ninnt\ we write $\pstr(n)$ for the finite sequence $p_0, \ldots, p_{n-1}$ of basic dots in $V$. Notice that $p_0\geqc\ldots\geqc p_{n-1}$, by definition of points. Write $\psstr(n)$ for the \precc-trail 
of $\pstr(n)$.  A finite sequence $a\iz a_0\succ\ldots\succ a_{n-1}$ of basic dots in $V$ is called a \deff{\precc-trail} from $a_0$ to $a_{n-1}$ of length $n$, or simply a trail from $a_0$ to $a_{n-1}$ in \Vprect. The empty sequence is the unique trail of length $0$, and denoted \maxdotpt. The countable set of trails in \Vprect\ is denoted \Vpath, notice that $\Vpath\iz\{\psstr(n)\midd\ninn, p\inn\Vcal\}$. 

\parr
Let $a\iz a_0, \ldots, a_{n-1}$ and  $b\iz b_0, \ldots, b_{m-1}$ be trails in \Vprect\ such that $a_{n-1}\succ b_0$, then we write $a\starr b$ for the concatenation $a_0,\ldots, a_{n-1}, b_0\ldots b_{m-1}$ which is again a trail and so in \Vpath. (Hereby $a\starr\maxdotp$ and $\maxdotp\starr a$ are defined to equal $a$.).
\parr
The basic dots of our trail space are the trails in \Vprect. For trails $a\iz a_0, \ldots,$ $a_{n-1}$ and  $b\iz b_0, \ldots, b_{m-1}$ 
we put: $a\leqcstar b$ iff there is a trail $c\inn\Vpath$ in such that $a\iz b\starr c$. We also put $a\aprtstar b$ iff $a_{n-1}\aprt b_{m-1}$. 
The natural space \Vpathnatt\ defined by the pre-natural space \Vpathpret\ is called the \deff{trail space}\ of \Vnatt.
\parr
Finally, a \leqc-morphism $f$ from \Vpathnatt\ to another natural space \Wnattot\ is called a \deff{trail morphism}\ (notation $\pthh$-morphism) from \Vnatt\ to \Wnattot. For a point $p\inn\Vcal$ we write $f(p)$ for the point of \Wcalt\ given by $f(\psstr(0)), f(\psstr(1)), \ldots$.
\edefi

\rem
From the pointwise perspective, one readily sees that \Vpathnatt\ is `just another representation' of \Vnatt. Differences in representation should be filtered out by the concept of `isomorphism'. This is the main reason for introducing trail morphisms, since \Vpathnatt\ is not always \leqc-isomorphic to \Vnatt\ (for an example consider the natural real numbers). In fact refinement morphisms preserve the lattice-order properties of the basic neighborhood system which is chosen for a specific representation. Due to the presence of an apartness/equivalence relation, these order properties are not always relevant since we can freely add or distract equivalent basic dots to our system with different lattice properties, without essentially changing the point space.
\erem

\thm
Let \Vnatt\ and \Vpathnatt, \Vpathpret\ be as in the above definition. Then 
\be
\item[(i)] \Vpathpret\ is a pre-natural space and \Vpathnatt\ is a natural space.
\item[(ii)]
\Vpathnatt\ is homeomorphic to \Vnatt\ as a topological space. A homeomorphism is induced by the \pthh-morphism \idptht\ from \Vnatt\ to \Vpathnatt\ given by $\idpth(p)\iz \psstr(0), \psstr(1), \psstr(2),\ldots\inn\Vcalpth$ for $p\inn\Vcal$ (as a refinement morphism \idptht\ is the identity on \Vpatht, with $\idpth(a)\iz a$ for $a\inn\Vpath$). Its inverse homeomorphism is induced by the \leqc-morphism \idstrt\ from \Vpathnatt\ to \Vnatt\ which is defined by putting $\idstr(\maxdotp)\iz\maxdot$ and $\idstr(a)\iz a_{n}$ for a trail $a\iz a_0, \ldots, a_{n}$ in \Vpatht.
\item[(iii)]
Let $f$ be a \pthh-morphism from \Vnatt\ to \Wnattot. Then $f$ is continuous.
\ee
\ethm

\prf
See\WaaNat, 1.1.4.
\eprf

If $f$ is a \leqc-morphism from \Vnatot\ to \Wnattot, then $f\circ\idstr$ is by definition a \pthh-morphism from \Vnatot\ to \Wnattot, which is clearly equivalent to $f$ on \Vcalt. Therefore we will consider each \leqc-morphism to be a \pthh-morphism as well.

\sbsc{How to unglue spraids (definition)}\label{unglueing}
We show that any spraid \Vnatt\ can be unglued. The idea is to turn to the subspread of the trail space \Vpathnatt\ which is formed by the \sucz-trails in \Vleqct\ (instead of looking at the tree \Vpatht\ of all trails).

\defi
Let \Vnatt\ be a spraid derived from \Vprenatt. The \deff{unglueing}\ of $(\Vcal,\Topaprt\!)$ is the spread \Vunglnatt\ derived from the pre-natural space \Vunglpre, where $\Vungl=\{a\iz a_0, \ldots a_{n-1}\inn\Vpath\midd \ninn\midd a\ \mbox{is a}\ \mbox{\sucz-trail}$ $\mbox{and}\ n\geqq 1\rightarrow\grdd(a_0)\iz 1\}$.
\edefi

We leave it to the reader to verify that unglueing a spraid \Vnatt\ amounts to adding, for each $a\inn V$, a finite number of copies of $a$ such that each \sucz-trail from \maxdott\ to $a$ is represented by one of the copies. These copies all have $\grdd(a)$ as length in $(\Vungl, \leqcstar)$.

\xam
We consider the important spread \sigRt. To turn this spraid into an isomorphic spread, we unglue. We look at the \sucz-trails\ in \sigRptht\ which (if not equal to the empty sequence $\maxdot^*$) start with a basic interval in $\suczz(\maxdotR)\iz\{[m, m\pluz 2]\midd m\inn\Z\}$. So put:

\parr

$\sigRungl\isdef\{a\iz a_0,\ldots a_{n-1}\inn \sigRpth\midd \ninn\midd a\ \mbox{is a}\ \mbox{\sucz-trail and}\ n\geqq 1\rightarrow\grdd(a_0)\iz 1 \}$

\parr

Then \sigRunglt\ has as maximal dot $\maxdot^*$, and an example of a basic dot in \sigRunglt\ is the sequence $[0, 2], [1,2]$, which has as unglued twin the basic dot $[1, 3], [1,2]$.
For simplicity, we also write \sigRunglt\ for the spread derived from the pre-natural space $(\sigRungl, \leqcstar, \aprtstar)$, which is the unglued version of \sigRt.
\exam

\sbsc{Proof of theorem \ref{metcompfan}}\label{prfmetcompfan}
For the proof of theorem \ref{metcompfan} we copy some of the representation theory in \WaaNat. 

\thm (from \ref{metcompfan})
Every compact metric space is homeomorphic to a fan.
\ethm

\prf
By the previous paragraph (\ref{unglueing}) it suffices to show that \xdt\ is homeomorphic to a fann (since this can be unglued to a fan).  Let \xdt\ be a compact metric space, meaning \xdt\ is totally bounded and complete. Using countable choice (\aczz) we can determine a sequence of finitely enumerable subsets $(C_i)_{i\in\N}$ of $X$ (for each $i$ with cardinality less than $k_i\pluz 2$, and given by $C_i\iz\{c_{i,j}\midd j\leqq{k_i}\}$) such that

\stArtabb
$\all x\inn X\,\all\ninn\dris c\inn C_n [d(x,c)\smlr 2^{-n-2}]$
\etab

Then $A\iz\{a_n\midd\ninn\}\iz\bigcupp_i C_i$ is a dense subset of \xdt. To define the required fann, we start with the (too large) set of basic dots $V=\{B(a_n, 2^{-s})\midd n,\sinn\}\cupp\{\maxdotV\}$. The technical trouble now is to define \aprt\ and \leqct\ constructively, since in general even for $s\bygr t$ the containment relation $B(a_n, 2^{-s})\subseteqq B(a_m,2^{-t})$ is not decidable. However, this containment relation has an enumerable subrelation which also does the trick. This because for all $(a_n, s)$ and $(a_m, t)$ with $s\bygr t$ there is $k\inn\{0,1\}$ such that:
\parr
$(k\iz 0 \weddge d(a_n, a_m)\smlr 2^{-t}\minuz 2^{-s})$ or $(k\iz 1\weddge d(a_n, a_m)\bygr 2^{-t}\minuz 2^{-s}\minuz 2^{-2s})$
\parr
Using \aczz\ (countable choice) we can define a function $h$ fulfilling the above statement. Now we put $B(a_n, 2^{-s})\precc B(a_m, 2^{-t})$ iff $h((a_n, s), (a_m, t))\iz 0$.
Likewise we define \aprt, since for all $(a_n, s)$ and $(a_m, t)$ there is $l\inn\{0,1\}$ such that:\parr
$(l\iz 0 \weddge d(a_n, a_m)\smlr 2^{-s\!}\pluz 2^{-t\!}\pluz 2^{-s-t} )$ or $(l\iz 1\weddge d(a_n, a_m)\bygr 2^{-s\!}\pluz 2^{-t\!}\pluz 2^{-s-t-1} )$
\parr
Using \aczz\ we can define a function $g$ fulfilling the above statement. Now we simply put $B(a_n, 2^{-s})\aprt B(a_m, 2^{-t})$ iff $g((a_n, s), (a_m, t))\iz 1$.
\parr
Taking as its basic dots the set $W=\{B(c_{i,j}, 2^{-i})\midd i\inn\N, j\leqq{k_i}\}\cupp\{\maxdotV\}$ we now form the required fann \Wnatt\ from the pre-natural space \Wprenatt. The verification that \Wnatt\ is a fann is relatively straightforward, by checking the conditions set in \stArt\ and the definitions of \leqct\ and \aprt. We define a homeomorphism $f$ from \Wnatt\ to \xdt\ as follows. For $w\iz w_0, w_1,\ldots \inn\Wcal$ with $w_0\notiz\maxdotV$ there is a unique Cauchy-sequence \enninnt\ in $A$ such that $e_n$ is the center of the metric ball formed by $w_n$. Now define $f(w)\isdef \lim_d\enninn$ (we will shortly show $f$ is a homeomorphism). Using $f$ it is easy to define the relevant metric $d'$ on \Wcalt\ such that $(\Wcal, d')$ is homeomorphic to \xdt: simply put $d'(w,z)\isdef d(f(w),f(z))$. That this makes $f$ a homeomorphism from $(\Wcal, d')$ to \xdt is ensured by \stArt\ above.
\parr
To show that $f$ is also a homeomorphism from \Wnatt\ to \xdt\ (thus finishing the proof) it therefore suffices to show that the $d'$-induced metric topology on \Wcalt\ coincides with the apartness topology. So let \Ucalt\ be open in \Wnatt, we need to show that \Ucalt\ is $d'$-open. For this let $x\inn\Ucal$. Consider $y\iz f(x)\inn X$. By our special constraints of \stArt\ above and the definition of \leqct, we can find a $w\inn\Wcal, w\equivv x$ such that for each \ninnt\ the $d$-ball $w_n$ is a neighborhood of $y$. Since \Ucalt\ is open and $w\equivv x$, we find \ninnt\ such that $\hattr{w_n}\subseteqq\Ucal$. From this it follows that for some \sinnt\ the $d'$-ball $B(x,2^{-s})$ is contained in \Ucalt.
\eprf

\end{document}